\documentclass{article}

\usepackage{amsfonts}
\usepackage{amssymb}
\usepackage{amsthm}
\usepackage{amsmath}
\usepackage{newlfont}
\usepackage{graphicx}

\usepackage{epsfig}

\usepackage{multicol}
\usepackage[T1]{fontenc}
\newtheorem {proposition}{Proposition}[section]
\newtheorem {theorem}{Theorem}[section]
\newtheorem {lemma}{Lemma}[section]
\newtheorem {example}{Example}[section]
\newtheorem {definition}{Definition}[section]

\newtheorem {corollary}{Corollary}[section]

\author{{\DJ{}or\dj{}e Barali\' {c}}\\ {\small Mathematical Institute SASA}\\[-2mm] {\small Belgrade, Serbia} \and Igor Spasojevi\'{c}\\ {\small Cambridge University}\\[-2mm] {\small Cambridge, United Kingdom}
\\[-2mm]}

\title{Illumination of Pascal's Hexagrammum and Octagrammum Mysticum}
\date{}
\begin{document}
\maketitle

\begin{abstract} We prove general results which include classical facts about 60 Pascal's lines as special cases.
Along similar lines we establish  analogous results about
configurations of 2520 conics arising from Mystic Octagon. We
offer a more combinatorial outlook on these results and their dual
statements. B\'{e}zout's theorem is the main tool, however its
application is guided by the empirical evidence and computer
experiments with program \textit{Cinderella}. We also emphasize a
connection with $k$-nets of algebraic curves.

\end{abstract}

\renewcommand{\thefootnote}{}
\footnotetext{This research was supported by the Grant 174020 of
the Ministry for Education and Science of the Republic of Serbia.}

\section{Introduction}

`Projective geometry is all geometry' was a dictum of 19th century
mathematics. Great masters Pascal, Steiner, Cayley, Salmon and
others have discovered beautiful theorems about interesting
geometric configurations. and some of them are well-known to wide
mathematical community today. There are many papers that still
treat this topics so it is not easy to set plot for our story.

Pascal, with his `Hexagrammum Misticum' (Mystic Hexagon)
\cite{Pasc}, already in 1639 found a necessary and sufficient
condition for six points to lay on the same conic and started
construction of the Hexagrammum Misticum (Mystic Hexagon)
Configuration. Steiner was the first who drew attention of
mathematicians to the complete figure obtained by joining in all
possible ways six points on a conic. There are 60 different ways
to do that so there exist 60 Pascal lines. Steiner \cite{Stei}
proved at the beginning of the XIX century that the 60 Pascal
lines are concurrent by triples in 20 points, known today as
Steiner points. Kirkman's \cite{Kirkman} main contribution was the
observation that Pascal lines meet also by triples over 60 points
`Kirkman's points' which are different from Steiner points and
form a $(60)_3$-type configuration.

Pl\"{u}cker showed that 20 Steiner points lie in fours on 15
lines, three through each point. This lines are called
Steiner-Pl\"{u}cker or just Steiner lines. Cayley and Salmon
discovered that Kirkman points lie in threes on 20 Cayley-Salmon
lines and that Cayley-Salmon lines meet in threes in 15 Salmon
points.

Veronese in a remarkable paper \cite{Vero} proved 'Veronese's
Decomposition Theorem' which states that $(60)_3$-type
configuration splits properly into six Desargues Configurations of
the type ${(10)_3}_1$'s. Veronese's proof relies on a clever
choice of straight lines involved and on a skilful use of
Desargues two triangle Theorem, an idea which goes back to
Kirkman. He also proved that there are infinitely many systems
consisting of sixty lines and points.

Following the classical works, many papers about the Pascal
configuration have been written. Some of them extended the theorem
to higher dimensions, see \cite{Bott} and \cite{Cour}. Other were
focused on finding easier and more elementary proofs. Since a lot
of lines and points appear in configurations, it is not so easy to
find clear notation which would guide and explain which hexagons
lead to interesting results. This was considerably clarified in
papers \cite{Ladd} and \cite{Clark} where all lines and points
appearing in configuration are connected with certain subgroups of
the permutation group $S_6$.

\textit{Octagrammum Mysticum} was originally appeared as a problem
in \cite{Wilk}. It has been studied in the last 140 years, but not
as much as the Hexagrammum Mysticum. Note that both results as
well as their duals admit \cite{Katz} a common generalization to
the case of a  $2n$-gon inscribed in a conic.

\medskip

Our objective is to establish some new results (Theorems
\ref{8con}, \ref{mysconin}, \ref{gslc}) about Hexagrammum Mysticum
and Octagrammum Mysticum. We also discuss the possibility of
connecting these results with more recent developments in
combinatorics and algebraic topology (Proposition \ref{moj1}).

In Section 2 we state some classical theorems about algebraic
curves.

In Section 3 we prove some general results about Pascal lines
using ideas of elementary algebraic geometry. From here it is easy
to deduce many interesting results involving lines, conics and
cubics passing through vertices of hexagons.

In Section 4 we continue in similar way, by proving that is
possible to produce many interesting Steiner lines. Also we are
proving a theorem about Salmon-Cayley line and discover some new
remarkable conics and cubics passing through points in
configurations.

In Section 5 we are studying octagons inscribed in conics. Our
attention is focused on conics which arise from permutations of
vertices of the octagon. We establish interesting facts about
these conics that are directly analogous to Steiner and Kirkman
points. Also we prove a result that generalizes the notion of the
Steiner line.

In Section 6 we state corresponding dual statements of previously
proved theorems.

In Section 7 we study the degenerate cases of some theorems
concerning the mystic hexagon and octagon.

In Section 8 we describe how these constructions could be used to
produce examples in the theory of arrangements and how we could
associate some combinatorial and algebraic objects to them.

In Section 9 we state and briefly discuss some related open
problems.

\section{Intersections of Algebraic Curves}

Our guiding principle is that a geometric problem often can be
interpreted as a question about intersections of carefully chosen
algebraic curves. This approach gives more flexibility and
provides easy proofs of geometrical facts involving mystic hexagon
and octagon.

Theory of algebraic curves is well understood and developed area
of mathematics. There are many monographes about this topic, for
example \cite{Grif}, \cite{gh} and \cite{Kir}. However, the
emphasis in our paper is on combinatorial constructions, motivated
by the experiments with program  \textit{Cinderella}, so all we
need in most of the constructions is a weak form of B\'{e}zout's
theorem and its immediate consequences, see \cite{Kir} Section
3.1. We also use the theorem of Cayley \cite{Cay} and Bacharach
\cite{Bach} in the form stated in \cite{Boro}.

\begin{theorem}[The Weak Form of B\'{e}zout's theorem]\label{wbt}
If two projective curves $\mathcal{C}$ and $\mathcal{D}$ in
$\mathbb{C} P^2$ of degrees $n$ and $m$ respectively have no
common component then $\left| C\cap D\right| \leq n\cdot m$
points.
\end{theorem}

\begin{corollary}\label{vazno} If two projective curves $\mathcal{C}$ and
$\mathcal{D}$ in $\mathbb{C} P^2$ of degree $n$ intersect at
exactly $n^2$ points and if $n \cdot m$ of these points lie on
irreducible curve $\mathcal{E}$ of degree $m<n$, then the
remaining $n \cdot (n-m)$ points lie on curve of degree at most
$n-m$.
\end{corollary}

\begin{theorem}[The Cayley-Bacharach Theorem]\label{cbt}Let $\mathcal{A}$ and $\mathcal{B}$ be two algebraic curves in $\mathbb{C}
P^2$ of degrees $p$ and $q$ respectively which intersect at
$p\cdot q$ distinct points. Let $E\subset\mathbb{C} P^2$ be
algebraic curve of degree $r\leq p+q-3$ passing through $p\, q-1$
points of $\mathcal{A}\cap\mathcal{B}$. Then $\mathcal{E}$ passes
also through the last point of intersection.
\end{theorem}

\section{Generalized Steiner-Kirkman Points}

In this section we study Pascal's Mystic Hexagon.

\begin{theorem}[The Pascal's Line for Cubics]\label{plc} Let $A B C D E F$ be
hexagon inscribed in a conic $\mathcal{C}$ and let $\mathcal{D}_1$
and $\mathcal{D}_2$  be distinct cubics that pass through $A$,
$B$, $C$, $D$, $E$ and $F$. Let $P$, $Q$ and $R$ be three other
points of intersection of $\mathcal{D}_1$ and $\mathcal{D}_2$.
Then the points $P$, $Q$ and $R$ are collinear.
\end{theorem}

\noindent {\bf Proof:} This is an immediate consequence of
corollary \ref{vazno}.\hfill $\square$

\medskip

The line arising in Theorem \ref{plc} will be referred to as the
\textit{generalized Pascal line}.

\begin{figure}[h!h!]
\centerline{\epsfig{figure=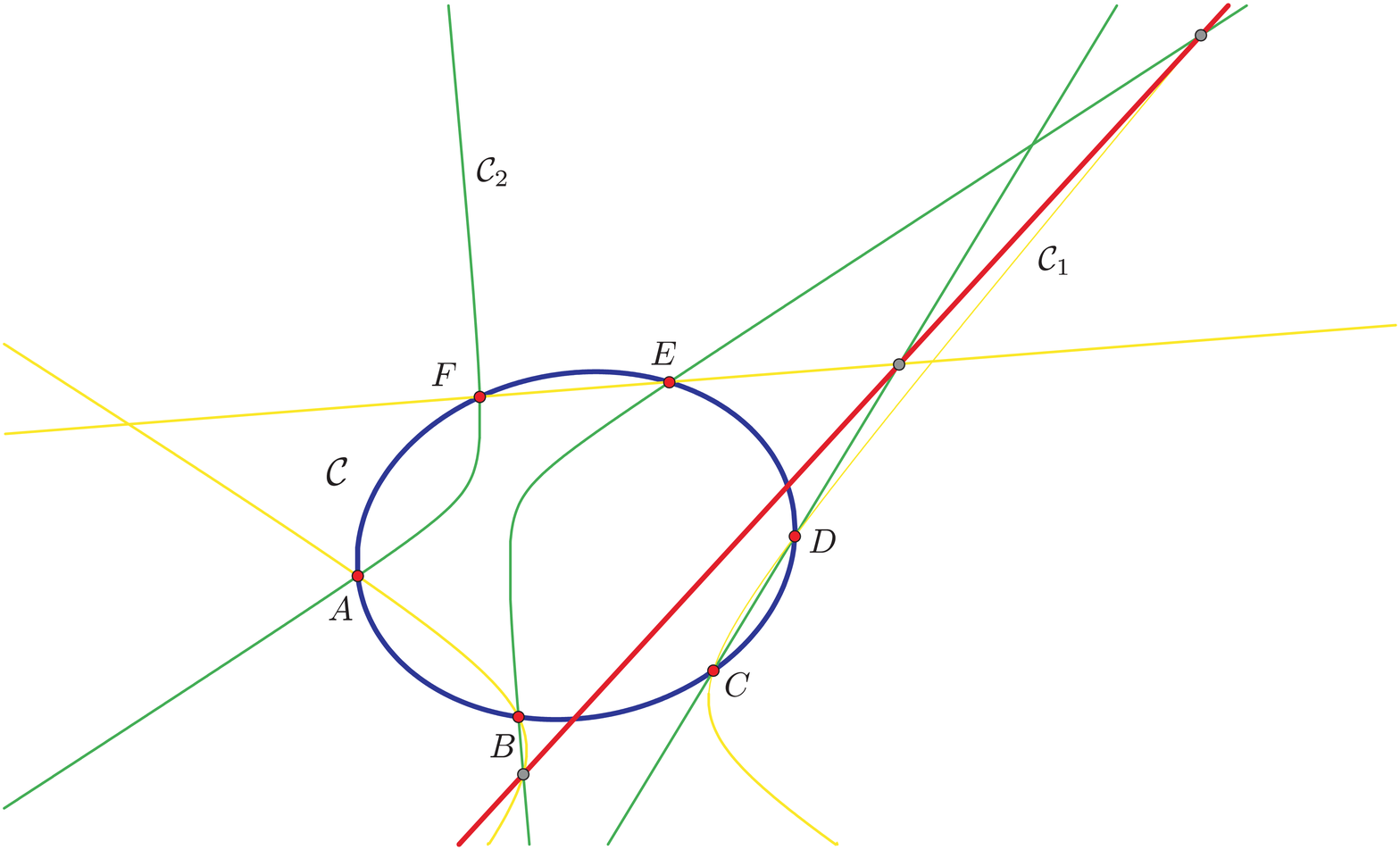, width=0.8\textwidth}}
\caption{Proposition \ref{slika1}} \label{prop2}
\end{figure}

From Theorem \ref{plc} easily follows:

\begin{proposition}\label{slika1} Let $A B C D E F$ be a
hexagon inscribed in a conic $\mathcal{C}$ and let $\mathcal{C}_1$
be a conic through points $A$, $B$, $C$ and $D$, $\mathcal{C}_2$ a
conic through points $A$, $B$, $E$ and $F$. Then the intersection
point of lines $C D$ and $E F$, and two intersection points of
conics $\mathcal{C}_1$ and $\mathcal{C}_2$ distinct from $A$ and
$B$ are collinear (see Figure \ref{prop2}).
\end{proposition}

From Theorem \ref{plc} we see that is possible to obtain many
generalized Pascal lines.

\begin{theorem}[The generalized Steiner-Kirkman Point]\label{GSKP} Let $A B C D E F$ be
a hexagon inscribed in a conic $\mathcal{C}$ and let
$\mathcal{D}_1$, $\mathcal{D}_2$ and $\mathcal{D}_3$ be distinct
cubics that pass through $A$, $B$, $C$, $D$, $E$ and $F$. Let
$p_1$ be the generalized Pascal line for cubics $\mathcal{D}_2$
and $\mathcal{D}_3$, $p_2$ for cubics $\mathcal{D}_1$ and
$\mathcal{D}_3$ and $p_3$ for cubics $\mathcal{D}_1$ and
$\mathcal{D}_2$. Then the lines $p_1$, $p_2$ and $p_3$ belongs to
the same pencil of lines.
\end{theorem}

\noindent {\bf Proof:} Consider the curves $\mathcal{D}_1\cdot
p_1$ and $\mathcal{D}_2\cdot p_2$. They intersect in $16$ points.
The curve $\mathcal{D}_3\cdot p_3$ passes through $15$ of them. By
Cayley-Bacharach Theorem it passes also through 16th point, the
intersection of lines $p_1$ and $p_2$. Obviously this point could
not belong to $\mathcal{D}_3$ so it belongs to the line
$p_3$.\hfill $\square$

Many geometrical results arise as consequences of Theorem
\ref{GSKP} if the cubics $\mathcal{D}_1$, $\mathcal{D}_2$ and
$\mathcal{D}_3$ are given a special meaning.

\begin{proposition}[The classical Steiner point] The classical Pascal's
lines of hexa-gons $A B E D C F$, $A D E F C B$ and $A F E B C D$
intersect in one Steiner point.
\end{proposition}

\noindent {\bf Proof:} Take $\mathcal{D}_1=l (A B)\cdot l (D
E)\cdot l(C F)$, $\mathcal{D}_2=l(B E)\cdot l(C D)\cdot l(A F)$
and $\mathcal{D}_3=l (A D)\cdot l (B C) \cdot l (E F)$.\hfill
$\square$

\begin{proposition}[The classical Kirkman point] The classical Pascal's
lines of hexagons $A B F D C E$, $A E F B D C$ and $A B D F E C$
intersect in one Kirkman point.
\end{proposition}

\noindent {\bf Proof:} Take $\mathcal{D}_1=l (A B)\cdot l (D
F)\cdot l(C E)$, $\mathcal{D}_2=l(A E)\cdot l(B F)\cdot l(C D)$
and $\mathcal{D}_3=l (A C)\cdot l (B D) \cdot l (E F)$.\hfill
$\square$

\begin{figure}[h!h!]
\centerline{\epsfig{figure=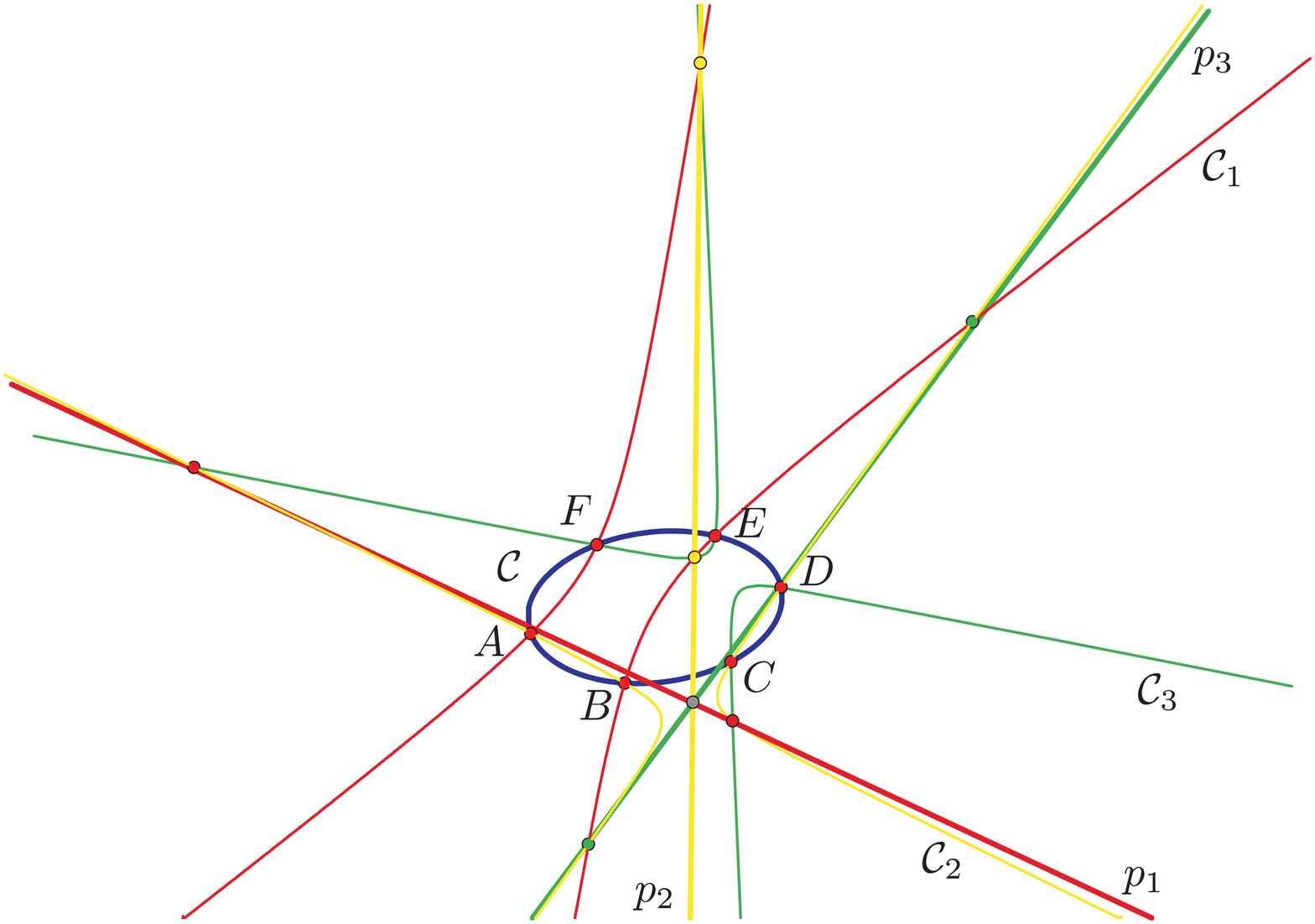, width=0.8\textwidth}}
\caption{Proposition \ref{slika2}} \label{prop1}
\end{figure}

\begin{proposition}\label{slika2} Let $A B C D E F$ be
a hexagon inscribed in a conic $\mathcal{C}$ and let
$\mathcal{C}_1$ be  a conic through points $A$, $B$, $C$ and $D$;
$\mathcal{C}_2$ a conic through points $A$, $B$, $E$ and $F$ and
$\mathcal{C}_3$ a conic through points $C$, $D$, $E$ and $F$. Let
$p_3$ be the line through the intersection points of
$\mathcal{C}_1$ and $\mathcal{C}_2$ distinct then $A$ and $B$ and
the lines $p_1$ and $p_2$ are defined in analogous way. Then the
lines $p_1$, $p_2$ and $p_3$ belong to the same pencil of lines
(see Figure \ref{prop1}).
\end{proposition}

\medskip

Theorem \ref{GSKP} allow us to define some interesting points
whose properties we discuss in the next section.

\begin{definition} \label{GSP} If in Theorem \ref{GSKP} we choose the cubic $\mathcal{D}_3$
arbitrary, $\mathcal{D}_1=l (A B)\cdot l (D F)\cdot l(C E)$ and
$\mathcal{D}_2=l (A D)\cdot l (B C) \cdot l (E F)$. Then obtained
point is called \textsl{generalized $\mathcal{D}_3$ Steiner
point}.
\end{definition}

We could take more restriction on the cubic $\mathcal{D}_3$. For
example to take for it product of conic and line and obtain the
results only involving conics and lines. However, our aim is to be
as general as it is possible for now.

\section{Generalized Steiner Lines, Salmon-Cayley Lines  and Salmon-Cayley Cubics}

In these sections we extend classical results about collinearities
among Steiner and Kirkman points to their generalized versions
introduced in Definition \ref{GSP}.

\begin{theorem}[The generalized Steiner line]\label{gsl} The four generalized $\mathcal{D}$ Steiner
points of hexagons $A B C D E F$, $A B D C E F$, $A D E B C F$ and
$A C E B D F$ lie on the same line.
\end{theorem}

\noindent {\bf Proof:} Let $p_1$ be the Pascal line for cubics
$\mathcal{D}$ and $l(A B)\cdot l(D F)\cdot l(C E)$; $p_2$ for
$\mathcal{D}$ and $l(A D)\cdot l(B C)\cdot l(E F)$; $p_3$ for
$\mathcal{D}$ and $l(A E)\cdot l(B D)\cdot l(C F)$ and $p_4$ for
$\mathcal{D}$ and $l(A C)\cdot l(B F)\cdot l(D E)$. Let $q_1$ be
the Pascal line for cubics $\mathcal{D}$ and $l(A C)\cdot l(E
F)\cdot l(B D)$; $q_2$ for $\mathcal{D}$ and $l(A B)\cdot l(C
F)\cdot l(D E)$; $q_3$ for $\mathcal{D}$ and $l(A D)\cdot l(B
F)\cdot l(C E)$ and $q_4$ for $\mathcal{D}$ and $l(A E)\cdot l(B
C)\cdot l(D F)$. The intersections $p_1\cap q_1$, $p_2\cap q_2$,
$p_3 \cap q_3$ and $p_4 \cap q_4$ are the associated generalized
Steiner points. Consider the curves $\mathcal{P}=p_1\cdot p_2\cdot
p_3\cdot p_4$ and $\mathcal{Q}=q_1\cdot q_2\cdot q_3\cdot q_4$.
Among $16$ intersection points of $\mathcal{P}\cap\mathcal{Q}$,
$12$ of them lie on the cubic $\mathcal{D}$ so by the Corollary
\ref{vazno}, the remaining four generalized Steiner points lie on
the same line. \hfill $\square$

From Theorem \ref{gsl} it is easy to obtain the classical results
about 20 Steiner points and 15 Steiner lines.

\begin{figure}[h!h!]
\centerline{\epsfig{figure=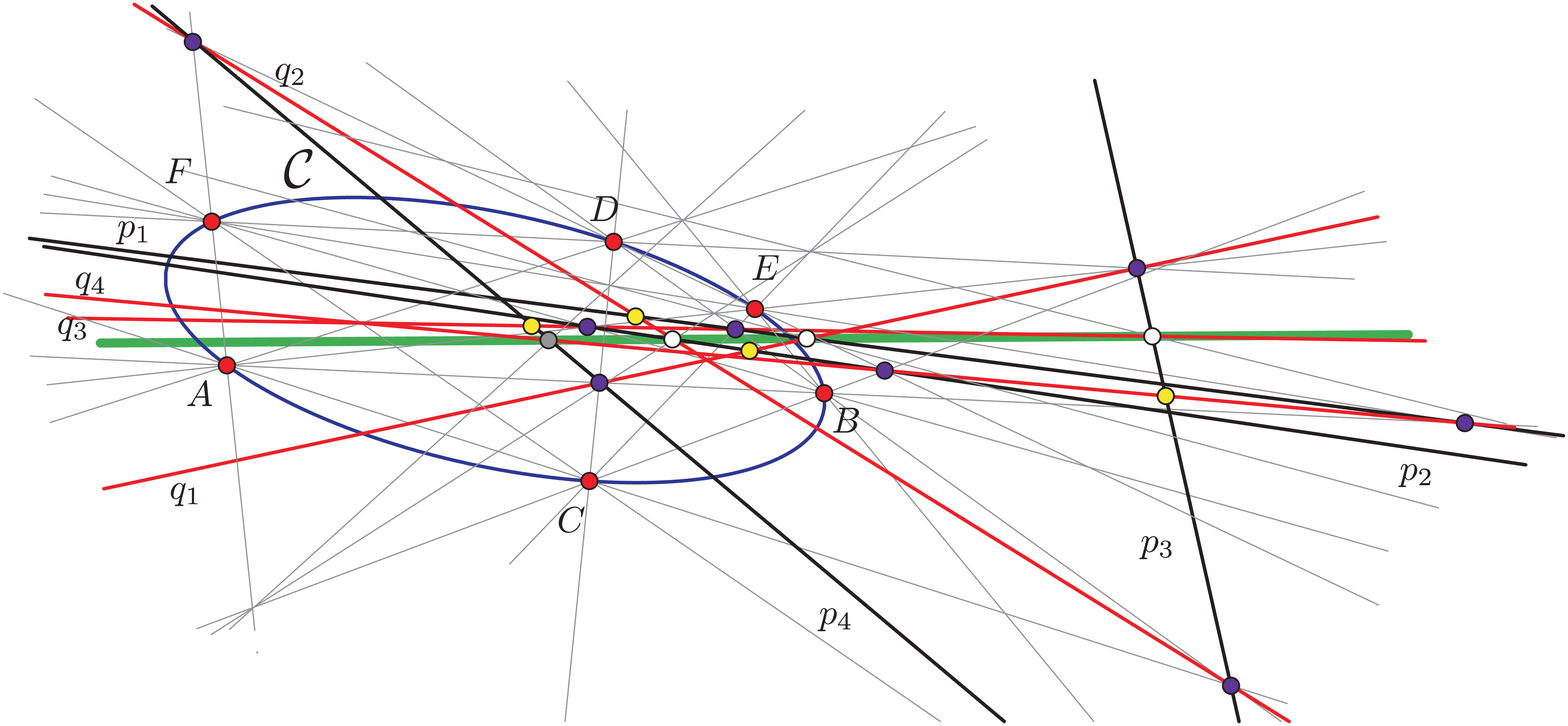, width=\textwidth}}
\caption{Theorem \ref{scl}} \label{salmoncayley}
\end{figure}

Now we proceed in the same manner to prove that 3 classical
Kirkman's points and the Steiner point lie on the same
Salmon-Cayley line (see Figure \ref{salmoncayley}). As we will
see, the proof of this fact is more complicated and we will
discover some interesting loci of points in Pascal configuration.

\begin{theorem}[The Salmon-Cayley Line]\label{scl} The Kirkman points of
hexagons \linebreak $A B D E F C$, $A C B F D E$ and $A C F D B E$
and the Steiner point of hexagon $A D C B E F$ lie on one line.
\end{theorem}

\noindent {\bf Proof:} Let $p_1$ be the Pascal line for cubics
$l(A B)\cdot l(D F)\cdot l(C E)$ and $l(A C)\cdot l(E F)\cdot l(B
D)$; $p_2$ for $l(A C)\cdot l(B F)\cdot l(D E)$ and $l(A E)\cdot
l(B C)\cdot l(D F)$; $p_3$ for $l(A C)\cdot l(B E)\cdot l(D F)$
and $l(A E)\cdot l(B D)\cdot l(C F)$ and $p_4$ for $l(A B)\cdot
l(D E)\cdot l(C F)$ and $l(A F)\cdot l(B E)\cdot l(C D)$. Let
$q_1$ be Pascal line for cubics $l(A B)\cdot l(D F)\cdot l(C E)$
and $l(A E)\cdot l(B F)\cdot l(C D)$; $q_2$ for $l(A C)\cdot l(B
F)\cdot l(D E)$ and $l(A F)\cdot l(B D)\cdot l(C E)$; $q_3$ for
$l(A C)\cdot l(B E)\cdot l(D F)$ and $l(A E)\cdot l(B D)\cdot l(C
F)$ and $q_4$ for $l(A B)\cdot l(D E)\cdot l(C F)$ and $l(A
D)\cdot l(B C)\cdot l(E F)$. The intersections $p_1\cap q_1$,
$p_2\cap q_2$, $p_3 \cap q_3$ and $p_4 \cap q_4$ are Kirkman and
Steiner points. Consider the curves $\mathcal{P}=p_1\cdot p_2\cdot
p_3\cdot p_4$ and $\mathcal{Q}=q_1\cdot q_2\cdot q_3\cdot q_4$.
The rest of the proof we split in two parts, Corollaries
\ref{2k1sp} and \ref{3kp}.

\begin{lemma}\label{lema2} The points $p_1\cap q_2$, $p_1\cap q_4$,
$p_2\cap q_1$, $p_2\cap q_4$, $p_4\cap q_1$ and $p_4\cap q_2$ lie
on the same conic.
\end{lemma}

\begin{figure}[h!h!]
\centerline{\epsfig{figure=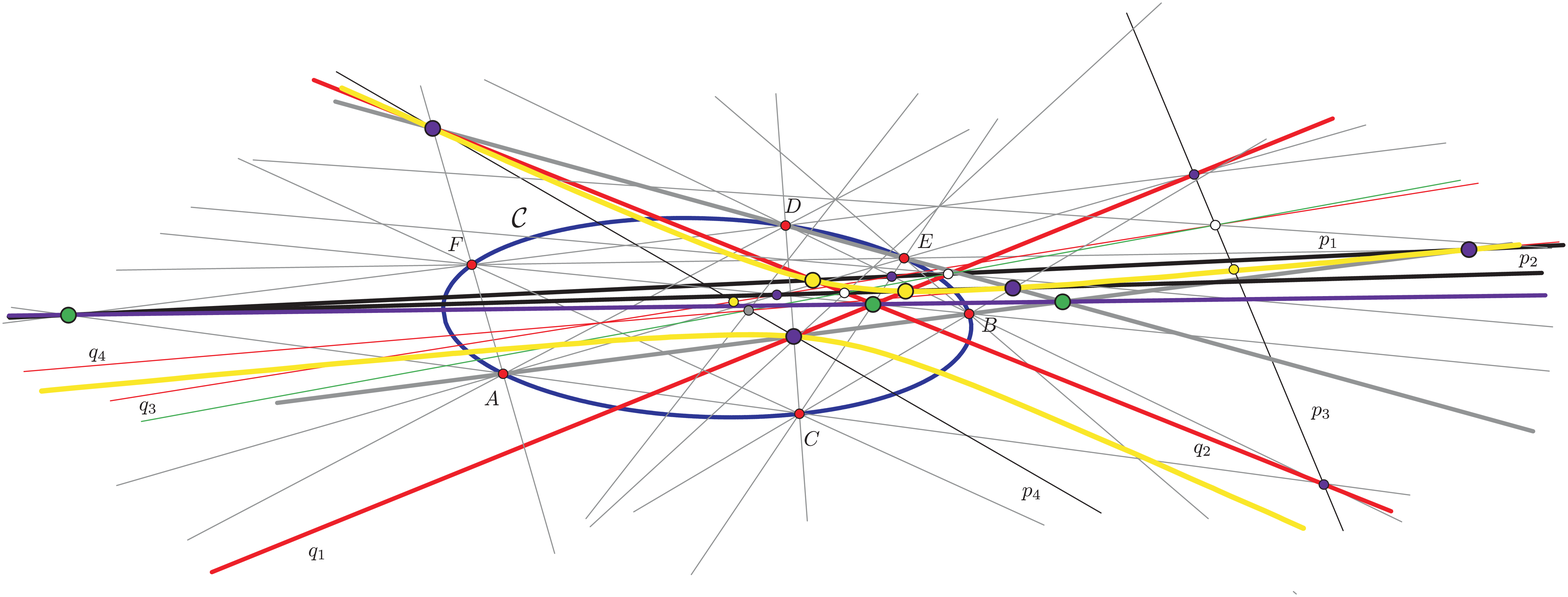, width=\textwidth}}
\caption{Lemma \ref{lema2}} \label{konikadokaz}
\end{figure}

\noindent {\bf Proof:} Look at the cubics $p_2 \cdot q_2\cdot A B$
and $p_1\cdot q_1\cdot D E$ (see Figure \ref{konikadokaz}). The
points $p_1\cap p_2$, $q_1\cap q_2$ and $A B\cap D E$ lie on the
same Pascal line and claim follows from the Corollary \ref{vazno}.
\hfill $\square$

\medskip

\begin{corollary} \label{2k1sp} Two Kirkman points $p_1\cap q_1$ and
$p_2\cap q_2$ and the Steiner point $p_4\cap q_4$ are collinear.
\end{corollary}

\medskip

\begin{lemma}\label{lema3} The points $p_3\cap q_1$, $p_1\cap q_3$,
$p_2\cap q_3$, $p_3\cap q_2$, $C$ and $F$ lie on the same conic.
\end{lemma}

\noindent {\bf Proof:} The points $C E\cap D F$, $p_2\cap q_2$ and
$A C\cap B F$ lie on the same Pascal line and the claim follows
from \ref{vazno}. \hfill $\square$

\begin{corollary} \label{wq} The points $A C \cap D F$,
$B F\cap C E$ and the intersection points of lines $l (p_3\cap
q_1, p_1\cap q_3)$ and $l (p_2\cap q_3, p_3\cap q_2)$ are
collinear.
\end{corollary}

\begin{lemma}\label{lema4} The points $p_3\cap q_1$, $p_1\cap q_3$,
$p_2\cap q_3$, $p_3\cap q_2$, $p_1\cap q_2$ and $p_2\cap q_1$ lie
on the same conic.
\end{lemma}

\noindent {\bf Proof:} Corollaries \ref{wq} and \ref{vazno} imply
the statement. \hfill $\square$

\begin{corollary} \label{3kp} Three Kirkman points $p_1\cap q_1$,
$p_2\cap q_2$ and $p_3\cap q_3$ are collinear.
\end{corollary}

Theorem about Salmon-Cayley line implies the following statement:

\begin{proposition}\label{propokubika} The points $p_1\cap q_2$, $p_1\cap q_3$, $p_1\cap
q_4$, $p_2\cap q_1$, $p_2\cap q_3$, $p_2\cap q_4$, $p_3\cap q_1$,
$p_3\cap q_2$, $p_3\cap q_4$, $p_4\cap q_1$, $p_4\cap q_2$ and
$p_4\cap q_3$ lie on one cubic.
\end{proposition}

This cubic we will be called \textit{Salmon-Cayley cubic} of
Pascal hexagon.

\begin{proposition}\label{propokonikd} The points $p_3\cap q_1$, $p_1\cap q_3$, $p_2\cap q_3$, $p_3\cap q_2$,
$C$ and $F$ lie on the same conic.
\end{proposition}

Conics in Lemmas \ref{lema2} and \ref{lema4} will be called
\textit{Steiner} and \textit{Kirkman} conic, respectively. They
have also one interesting property.

\begin{proposition} Two intersection points of Steiner and Kirkman
conics distinct from $p_1\cap q_2$ and $q_1\cap p_2$ lie on the
line $C F$.
\end{proposition}

\section{Conics in Octagrammum Mysticum}

In this section we will study a configuration obtained by 8 points
inscribed in a conic. Recall that that a quartic is an algebraic
curve of degree 4.

\begin{theorem}[The general Octagrammum Mysticum]\label{8con} Let $A B C D E F G H$ be
an octagon inscribed in a conic $\mathcal{C}$ and let
$\mathcal{Q}_1$ and $\mathcal{Q}_2$ be distinct quartics that pass
through the points $A$, $B$, $C$, $D$, $E$, $F$, $G$  and $H$. Let
$L$, $M$, $N$, $O$, $P$, $Q$, $R$ and $S$ be eight other points of
the intersection of $\mathcal{Q}_1$ and $\mathcal{Q}_2$. Then
these eight points lie on the same conic.
\end{theorem}

\noindent {\bf Proof:} Analogous to the proof \ref{plc}. \hfill
$\square$

\begin{figure}[h!h!]
\centerline{\epsfig{figure=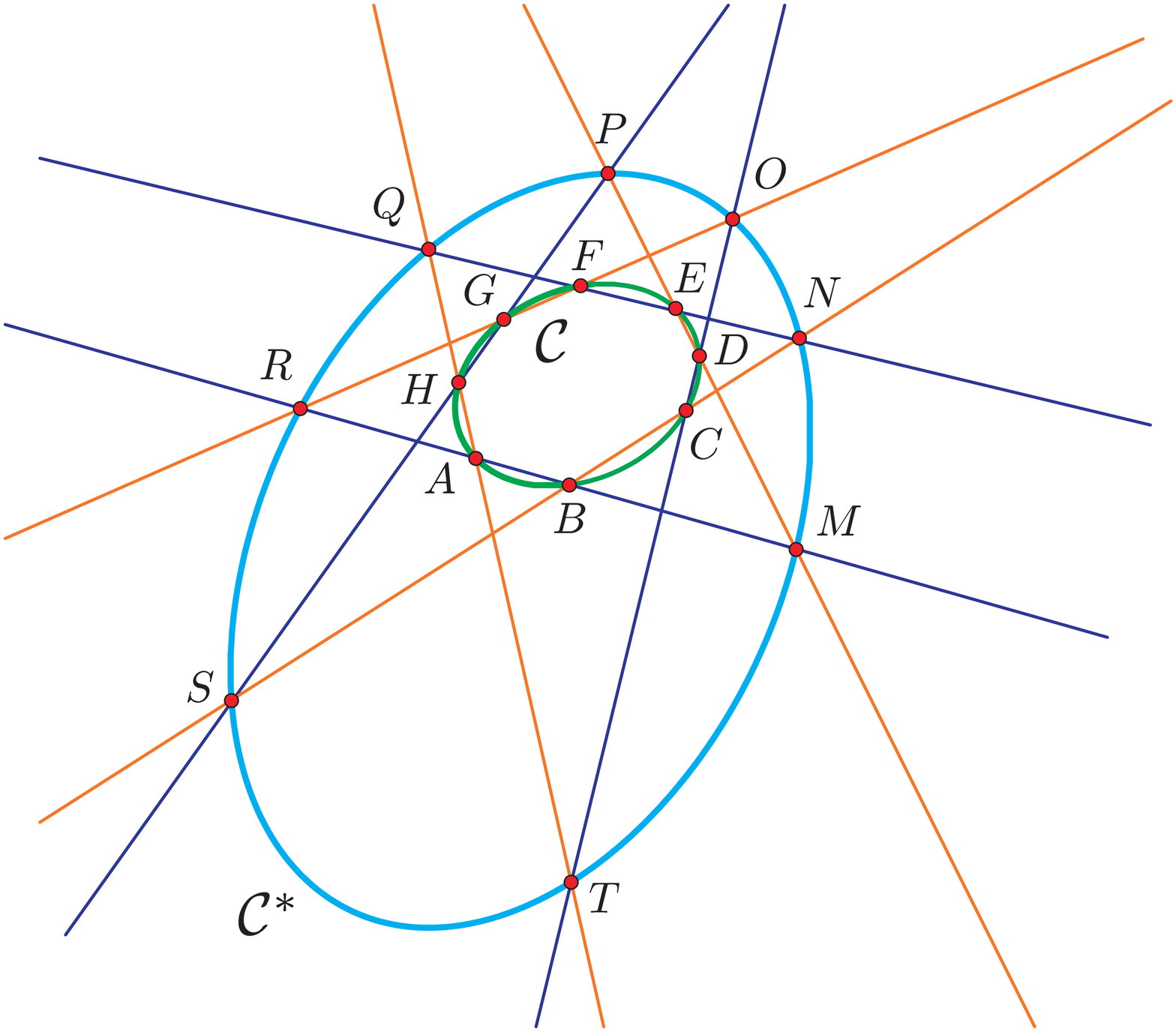,
width=0.5\textwidth}} \caption{Proposition \ref{8conc}}
\label{osmougaoelipsa}
\end{figure}

\begin{proposition}[The classical Octagrammum Mysticum]\label{8conc} Let $A B C D E F G H$ be
an octagon inscribed in a conic $\mathcal{C}$ and let the lines $A
B$, $C D$, $E F$ and $G H$ intersect the lines $B C$, $D E$, $F G$
and $H A$ in the points $M$, $N$, $O$, $P$, $Q$, $R$, $S$ and $T$.
Then the eight points $M$, $N$, $O$, $P$, $Q$, $R$, $S$ and $T$
lie on the same conic (see Figure \ref{osmougaoelipsa}).
\end{proposition}

Proposition \ref{8conc} suggests that is naturally to investigate
all possible ways to join $8$ points on a conic and corresponding
conics. However, there exist $\frac{8!}{8\cdot 2}=2520$ conics in
the classical configuration because reversing order and cyclic
permuting of vertices yields the same joining. Two conics in
$\mathbb{C}P^2$ intersects generally in $4$ points. The following
result says that some of these conics belong to the same pencil.

\begin{theorem}\label{mysconin} Let $A B C D E F G H$ be
an octagon inscribed in a conic $\mathcal{C}$ and assume that
$\mathcal{Q}_1$, $\mathcal{Q}_2$ and $\mathcal{Q}_3$ are distinct
quartics that pass through $A$, $B$, $C$, $D$, $E$, $F$, $G$ and
$H$. Let $\mathcal{C}_1$ be the mystic conic for quartics
$\mathcal{Q}_2$ and $\mathcal{Q}_3$, $\mathcal{C}_2$ for quartics
$\mathcal{Q}_1$ and $\mathcal{Q}_3$ and $\mathcal{C}_3$ for
quartics $\mathcal{Q}_1$ and $\mathcal{Q}_2$. Then the conics
$\mathcal{C}_1$, $\mathcal{C}_2$ and $\mathcal{C}_3$ belongs to
the same pencil of conics.
\end{theorem}

\noindent {\bf Proof:} Consider the curves $\mathcal{Q}_1\cdot
\mathcal{C}_1$ and $\mathcal{Q}_2\cdot \mathcal{C}_2$. They
intersect at $36$ points. The curve $\mathcal{Q}_3$ passes through
$24$ of them, $8$ vertices of octagon, $8$ points defining
$\mathcal{C}_1$ and $8$ points defining $\mathcal{C}_2$. By
Corollary \ref{vazno} the remaining $12$ points lie on degree $2$
curve which is obviously $\mathcal{C}_3$. \hfill $\square$

\begin{figure}[h!h!]
\centerline{\epsfig{figure=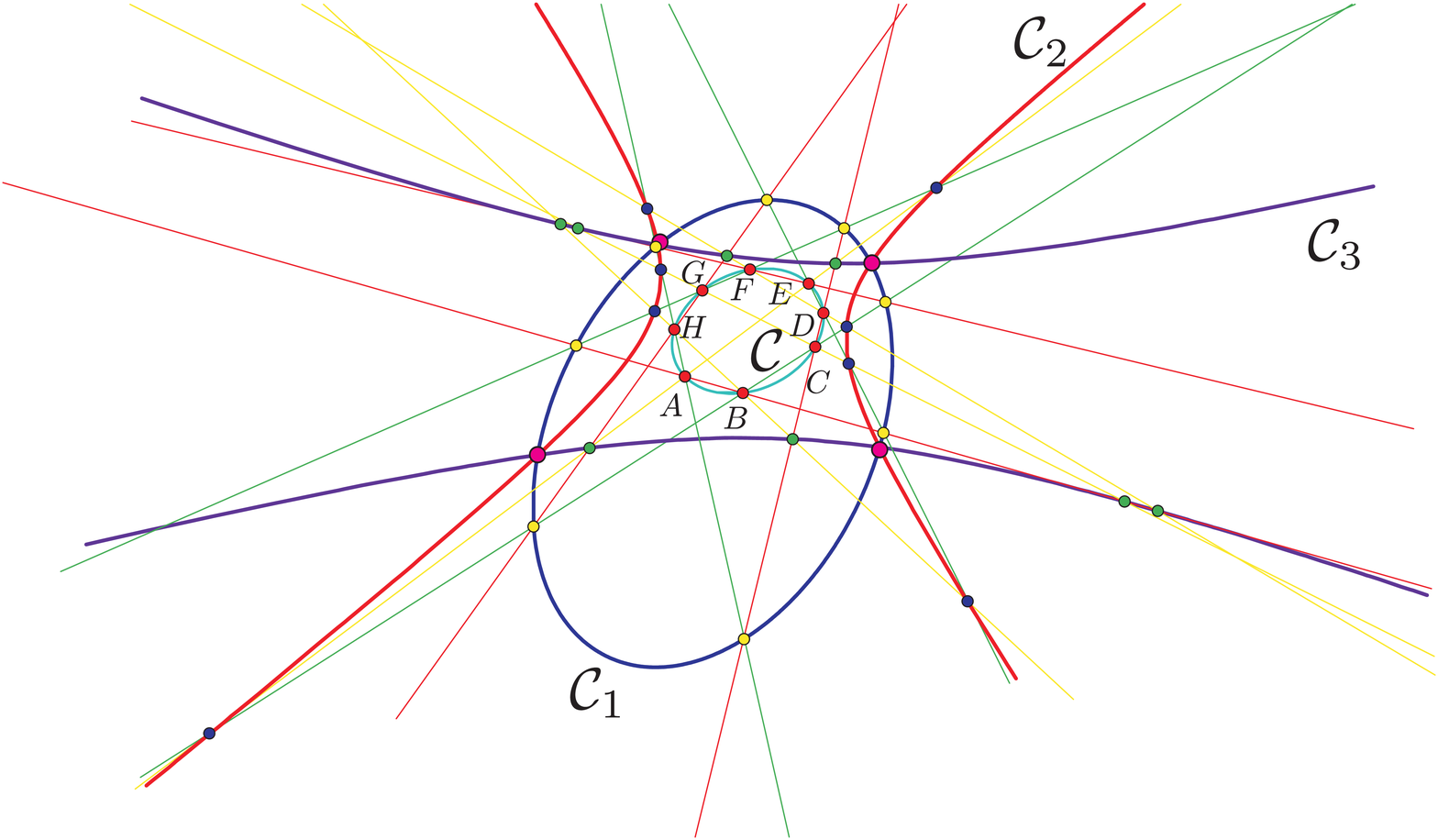, width=\textwidth}}
\caption{Theorem \ref{mysconin}} \label{konikeseseku}
\end{figure}

\medskip

Proposition \ref{2q} is a special case of Theorem \ref{8con}
corresponding to the case of two quadrilaterals inscribed in a
conic. There are $630$ such conics in the classical case of
quartics formed by 4 lines. We will see in Theorem \ref{gslc} why
this is interesting case.

\begin{proposition}\label{2q} Let $A B C D$ and  $E F G H$ be
quadrilaterals inscribed in a conic $\mathcal{C}$ and let the
lines $A B$, $C D$, $E F$ and $G H$ intersect the lines $B C$, $A
D$, $F G$ and $E H$ in the points $M$, $N$, $O$, $P$, $Q$, $R$,
$S$ and $T$. The eight points $M$, $N$, $O$, $P$, $Q$, $R$, $S$
and $T$ lie on the same conic.
\end{proposition}
 If we are interested only in a classical case, according to
Theorem \ref{mysconin}  there exist $28560$ pencils of conics such
that in each pencil lie $3$ of $2520$ conics, and each of $2520$
belongs to $34$ pencil of conics.

\begin{figure}[h!h!h!]
\centerline{\epsfig{figure=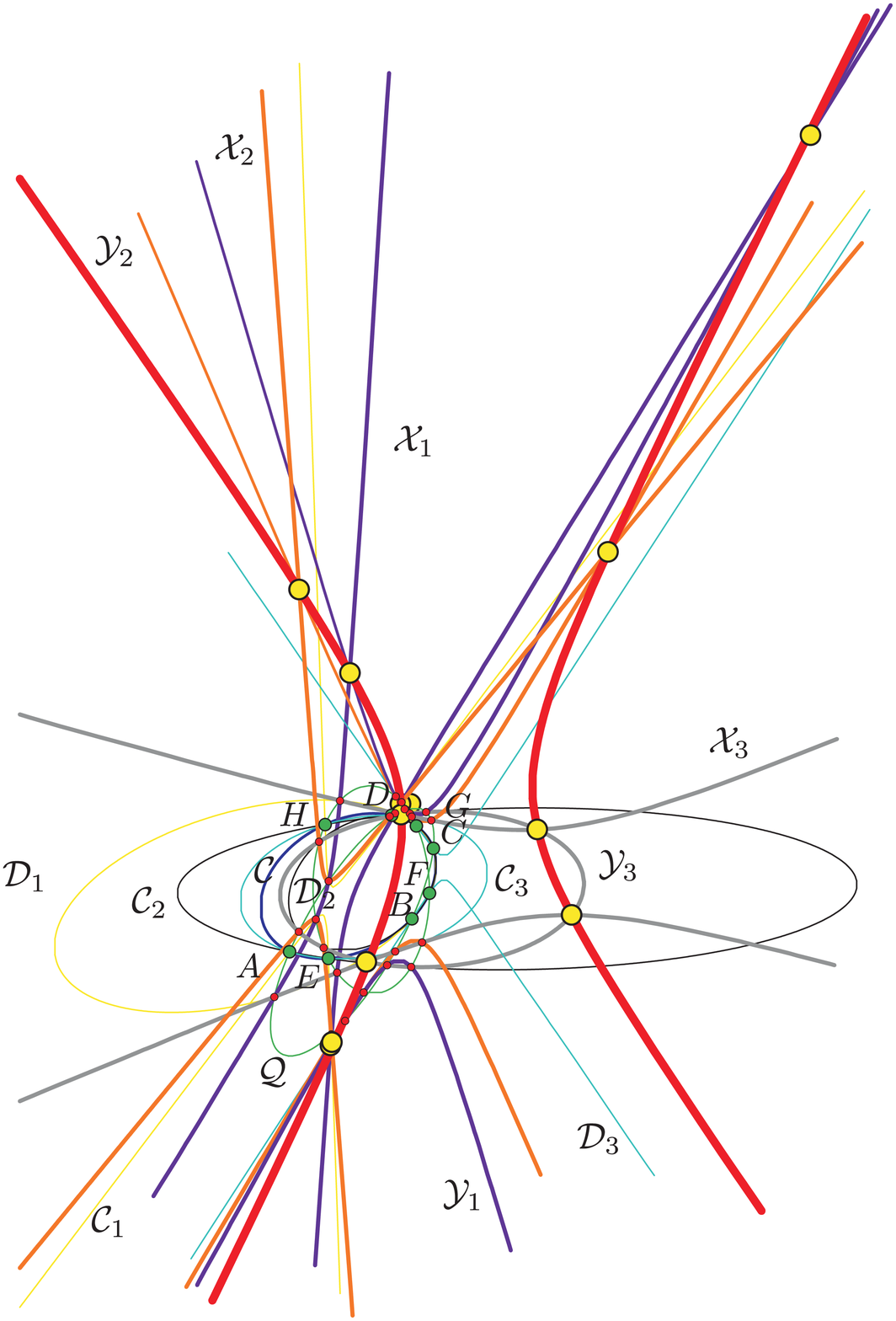, width=0.9\textwidth}}
\caption{Theorem \ref{gslc}} \label{stajnerkonika}
\end{figure}

The following theorem states that $3$ certain pencils have a
common conic. This result is analogous to Theorem \ref{gsl}.

\begin{theorem}\label{gslc} Let $\mathcal{Q}$ be a quartic passing
through $8$ vertices of mystic octagon, let $\mathcal{C}_1$,
$\mathcal{C}_2$ and $\mathcal{C}_3$ be three distinct conics
through points $A$, $B$, $C$ and $D$ and let $\mathcal{D}_1$,
$\mathcal{D}_2$ and $\mathcal{D}_3$ be three distinct conics
through points $E$, $F$, $G$ and $H$. Let $\mathcal{X}_1$ be the
mystic conic arising from the  curves $\mathcal{Q}$ and
$\mathcal{C}_1\cdot\mathcal{D}_1$ and $\mathcal{Y}_1$ be the
mystic conic arising from curves $\mathcal{Q}$ and
$\mathcal{C}_3\cdot\mathcal{D}_2$. The conics $\mathcal{X}_2$,
$\mathcal{X}_3$, $\mathcal{Y}_2$ and $\mathcal{Y}_3$ are defined
in analogous way. Then 12 intersection points of
$\mathcal{X}_1\cap\mathcal{Y}_1$, $\mathcal{X}_2\cap\mathcal{Y}_2$
and $\mathcal{X}_3\cap\mathcal{Y}_3$ lie on the same conic.
\end{theorem}

\noindent {\bf Proof:} Look at the curves $\mathcal{X}_1\cdot
\mathcal{X}_2\cdot \mathcal{X}_3$ and $\mathcal{Y}_1\cdot
\mathcal{Y}_2\cdot \mathcal{Y}_3$ (see Figure
\ref{stajnerkonika}). The quartic $\mathcal{Q}$ passes through
$24$ intersection points of this two curves so the remaining 12
must lie on the conic. \hfill $\square$

Theorem \ref{gslc} has many special cases because we are free to
choose special quartics and conics.

To conclude this section, we note that Theorems \ref{8con} and
\ref{mysconin} extend to the case of $2n$-gon inscribed in a conic
like in \cite{Katz}. The proof is the same as in the case of
octagon so we omit it.

\begin{theorem} Let $A_1 A_2\dots A_{2n}$ be
a $2n$-gon inscribed in a conic $\mathcal{C}$ and let
$\mathcal{Q}_1$ and $\mathcal{Q}_2$ be distinct degree $n$ curves
that pass through the vertices of $2n$-gon. Then the remaining
$n^2-2 n$ intersection points of $\mathcal{Q}_1\cap\mathcal{Q}_2$
lie on the curve of degree at most $n-2$.
\end{theorem}

\begin{theorem} Let $A_1 A_2\dots A_{2n}$ be a
$2n$-gon inscribed in a conic $\mathcal{C}$ and let
$\mathcal{Q}_1$, $\mathcal{Q}_2$ and $\mathcal{Q}_3$ be distinct
degree $n$ curves that pass through the vertices of $2n$-gon. Let
$\mathcal{C}_1$ be the mystic degree $n-2$ curve for
$\mathcal{Q}_2$ and $\mathcal{Q}_3$, $\mathcal{C}_2$ for
$\mathcal{Q}_1$ and $\mathcal{Q}_3$ and $\mathcal{C}_3$ for
$\mathcal{Q}_1$ and $\mathcal{Q}_2$. Then the curves
$\mathcal{C}_1$, $\mathcal{C}_2$ and $\mathcal{C}_3$ belong to the
same pencil of conics.
\end{theorem}

\section{Duality and Corresponding Results}

Duality between the points and the lines in projective geometry
allow us to formulate the corresponding dual theorems for conics
inscribed in a hexagon and an octagon. In this section we give
some interesting statements that are dual to the previously proved
theorems.

\begin{proposition} Let  $\mathcal{C}$ be a conic inscribed in a hexagon $A B C D E F$  and let
$\mathcal{C}_1$ be a conic touching the lines $A B$, $B C$, $C D$
and $A F$ and $\mathcal{C}_2$ a conic touching  the lines $A B$,
$D E$, $E F$ and $A F$. Then other two common tangents of
$\mathcal{C}_1$ and $\mathcal{C}_2$ and the line $C E$ intersect
at one point.
\end{proposition}

\begin{theorem}\label{presek} Let  $\mathcal{C}$ be a conic inscribed in a hexagon $A B C D E F$  and let
$\mathcal{C}_1$ be a conic touching the lines $A B$, $B C$, $C D$
and $D E $; $\mathcal{C}_2$ a conic  touching the lines $A B$, $B
C$, $E F$ and $A F$; and $\mathcal{C}_3$ a conic touching the
lines $C D$, $D E$, $E F$ and $A F$. The common tangents of
$\mathcal{C}_1$ and $\mathcal{C}_2$ distinct then $A B$ and $B C$
intersect at the point $P_3$. The points $P_1$ and $P_2$ are
defined in analogous way. Then the points $P_1$, $P_2$ and $P_3$
are collinear.
\end{theorem}

Dual statements for the classical Steiner and Kirkman points are
already known, see \cite{Ladd}. Duality argument could also be
applied in the case of octagon.

\begin{figure}[h!h!]
\centerline{\epsfig{figure=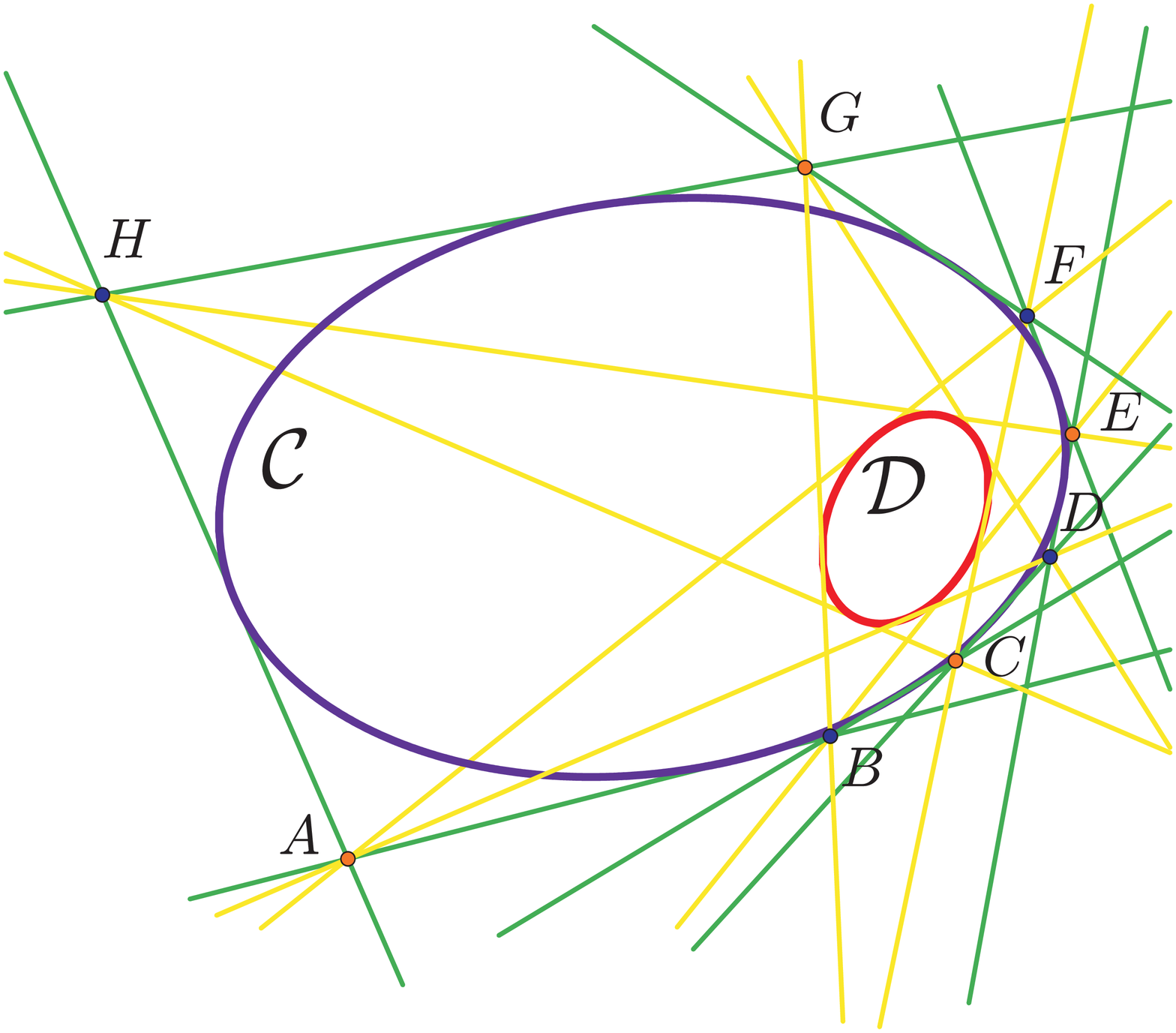, width=0.6\textwidth}}
\caption{Theorem \ref{os1}} \label{brijanosm}
\end{figure}

\begin{theorem}\label{os1} Let $\mathcal{C}$ be a conic inscribed in an octagon $A B C D E F
G H$. Let $\mathcal{M}$ be any octagon on the same vertices, then
there exist conic which tangents the sides of $\mathcal{M}$ (see
Figure \ref{brijanosm}).
\end{theorem}

\begin{figure}[h!h!]
\centerline{\epsfig{figure=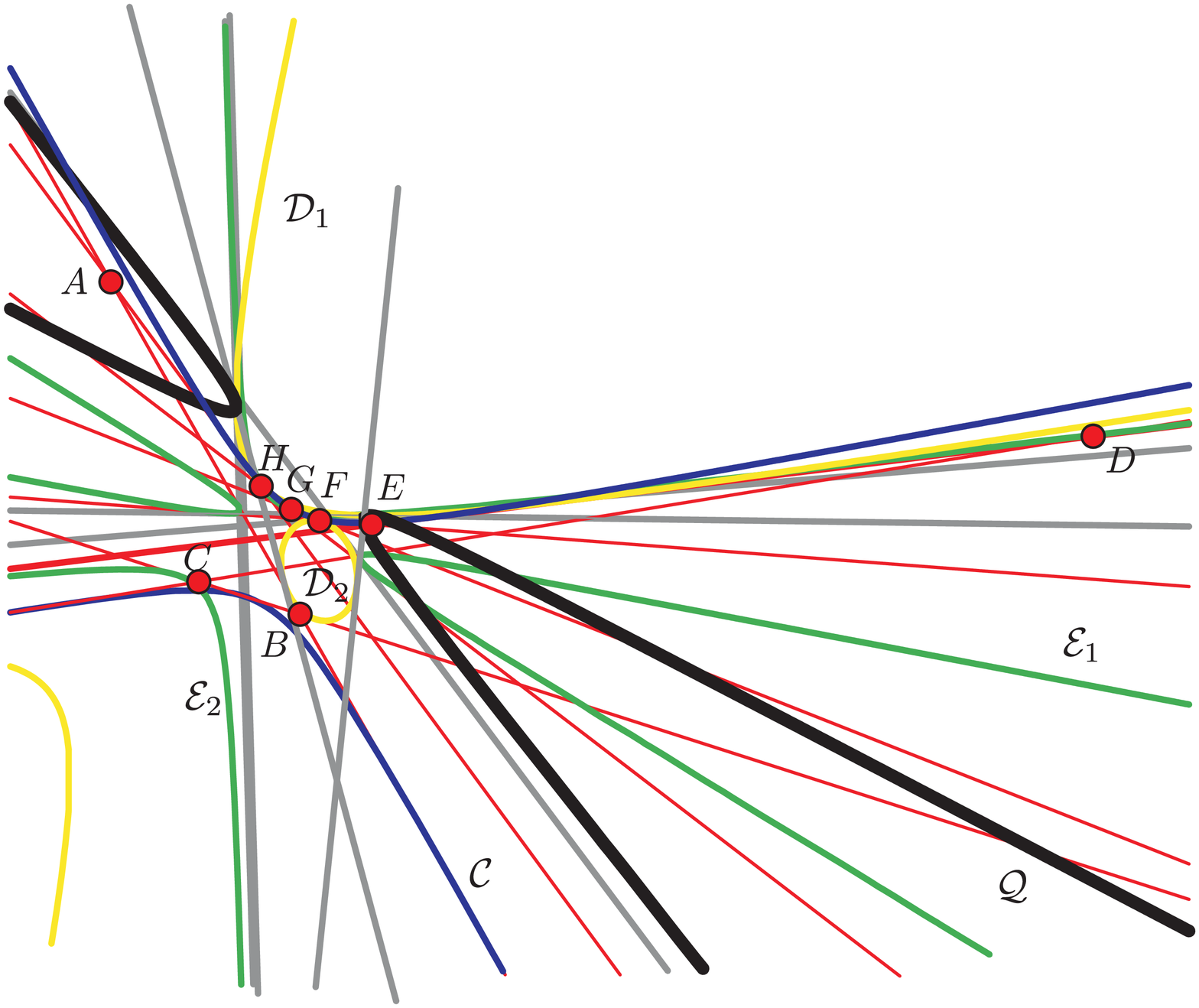,
width=0.8\textwidth}} \caption{Theorem \ref{lepo}}
\label{brijanosm1}
\end{figure}

\begin{theorem}\label{lepo} Let $\mathcal{C}$ be a conic inscribed in an octagon $A B C D E F
G H$. Let $\mathcal{D}_1$, $\mathcal{D}_2$ be conics touching some
four sides of the octagon and $\mathcal{E}_1$, $\mathcal{E}_2$ be
conics touching the remaining four sides of the octagon,
respectively. Then exist conic that touches the remaining $8$
common tangents of $\mathcal{D}_1$ and $\mathcal{E}_1$,
$\mathcal{D}_1$ and $\mathcal{E}_2$, $\mathcal{D}_2$ and
$\mathcal{E}_1$ and $\mathcal{D}_2$ and $\mathcal{E}_2$ (see
Figure \ref{brijanosm1}).
\end{theorem}

\begin{figure}[h!h!]
\centerline{\epsfig{figure=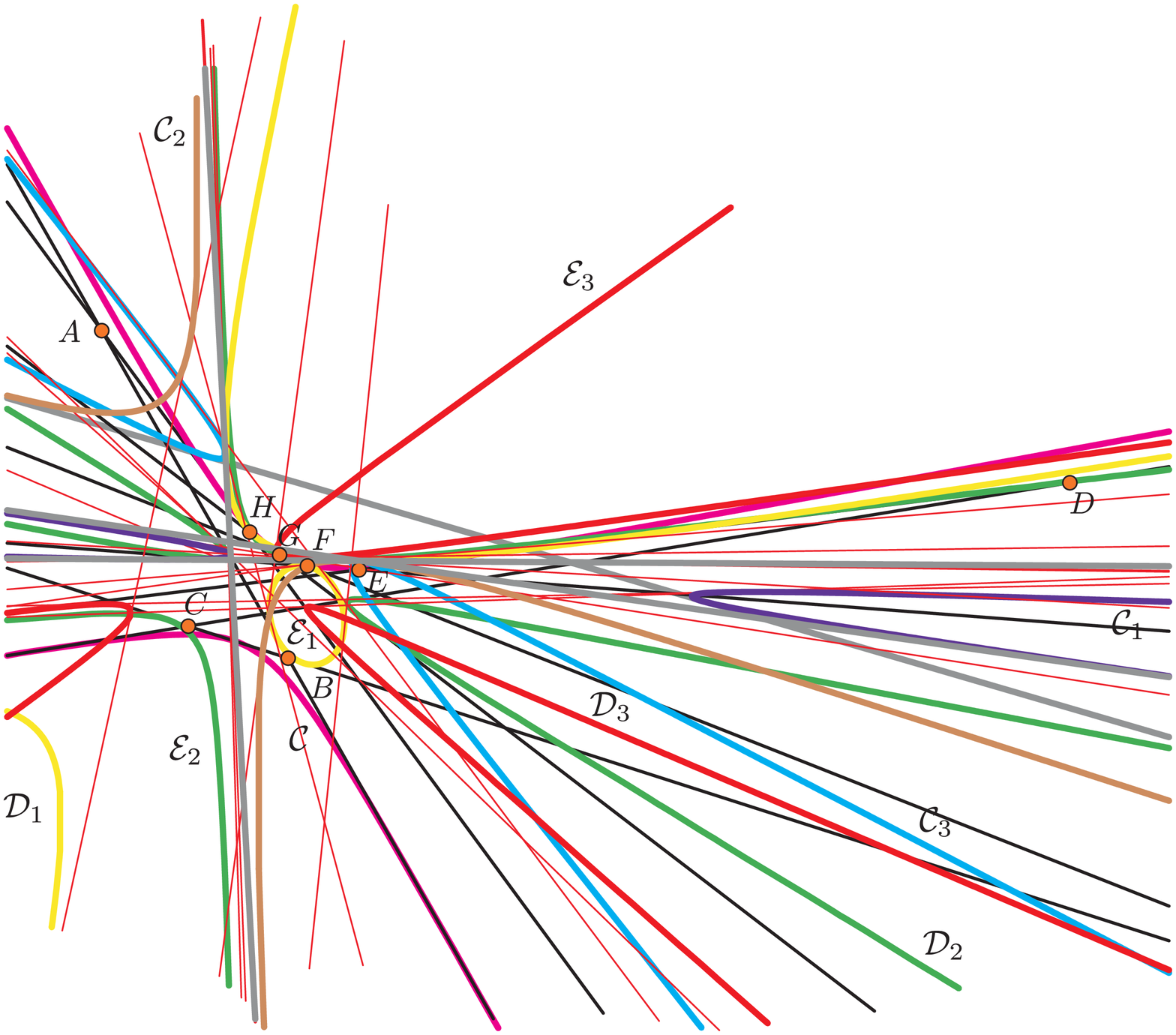,
width=\textwidth}} \caption{Theorem \ref{b1}} \label{brijanosm2}
\end{figure}

\begin{theorem}\label{b1} Let $\mathcal{C}$ be a conic inscribed in an octagon $A B C D E F
G H$. Let $\mathcal{D}_1$, $\mathcal{D}_2$, $\mathcal{D}_3$ be the
conics touching some four sides of the octagon each and
$\mathcal{E}_1$, $\mathcal{E}_2$, $\mathcal{E}_3$  be conics
touching the remaining four sides of the octagon, respectively.
Let $\mathcal{C}_1$ be the conic that touches remaining $8$ common
tangents of $\mathcal{D}_1$ and $\mathcal{E}_1$, $\mathcal{D}_1$
and $\mathcal{E}_2$, $\mathcal{D}_2$ and $\mathcal{E}_1$ and
$\mathcal{D}_2$ and $\mathcal{E}_2$. The conics $\mathcal{C}_2$
and $\mathcal{C}_3$ are defined in analogous way. Then there exist
the four lines each tangents $\mathcal{C}_1$, $\mathcal{C}_2$ and
$\mathcal{C}_3$ (see Figure \ref{brijanosm2}).
\end{theorem}

\begin{theorem}\label{b2} Let $\mathcal{C}$ be a conic inscribed in an octagon $A B C D E F
G H$, let $\textit{E}$ be a conic touching some four sides of the
octagon, and $\textit{F}$ be a conic touching the remaining four
sides of octagons. Let $\mathcal{C}_1$, $\mathcal{C}_2$ and
$\mathcal{C}_3$ be three distinct conics through the points $A B$,
$B C$, $C D$ and $D E$ and let $\mathcal{D}_1$, $\mathcal{D}_2$
and $\mathcal{D}_3$ be three distinct conics through the points $E
F$, $F G$, $G H$ and $H A$. Let $\mathcal{X}_1$ be the conic
arising from  Theorem \ref{b1} for   $\mathcal{E}$, $\mathcal{F}$,
$\mathcal{C}_1$ and $\mathcal{D}_1$, and let $\mathcal{Y}_1$ be
the conic arising from Theorem \ref{b1} for $\mathcal{E}$,
$\mathcal{F}$, $\mathcal{C}_3$ and $\mathcal{D}_2$. The conics
$\mathcal{X}_2$, $\mathcal{X}_3$, $\mathcal{Y}_2$ and
$\mathcal{Y}_3$ are defined in analogous way. Then there is conic
that tangents $12$ lines that are common tangents of
$\mathcal{X}_1$ and $\mathcal{Y}_1$; $\mathcal{X}_2$ and
$\mathcal{Y}_2$ and $\mathcal{X}_3$ and $\mathcal{Y}_3$.
\end{theorem}

\section{Degeneracy cases}

In the theorems about mystic hexagon and octagon is possible to
take the limit case when some vertices tends to some other
vertices. In that case configuration degenerates and we get
statements where both curves share a common tangent at that
vertex. In fact, these statements are special cases of the
previously proved results.

\begin{proposition} Let $A B C D E$ be a
pentagon inscribed in a conic $\mathcal{C}$ and let
$\mathcal{D}_1$ and $\mathcal{D}_2$ be distinct cubics that pass
through $A$, $B$, $C$, $D$, and $E$, such that there is common
tangent of $\mathcal{D}_1$ and $\mathcal{D}_2$ at $A$. Let $P$,
$Q$ and $R$ be three other points of intersections of
$\mathcal{D}_1$ and $\mathcal{D}_2$. Then the points $P$, $Q$ and
$R$ are collinear.
\end{proposition}

\begin{proposition} Let $A B C D E F G$ be
$7$-gon inscribed in a conic $\mathcal{C}$ and let $\mathcal{Q}_1$
and $\mathcal{Q}_2$ be distinct quartics that pass through $A$,
$B$, $C$, $D$, $E$, $F$ and $G$ such that there is common tangent
of $\mathcal{Q}_1$ and $\mathcal{Q}_2$ at $A$. Let $P$, $Q$, $R$,
$S$, $T$, $U$, $V$ and $W$ be 8 other points of the intersection
of $\mathcal{Q}_1$ and $\mathcal{Q}_2$. Then these points lie on
the sameconic.
\end{proposition}

\begin{figure}[h!h!]
\centerline{\epsfig{figure=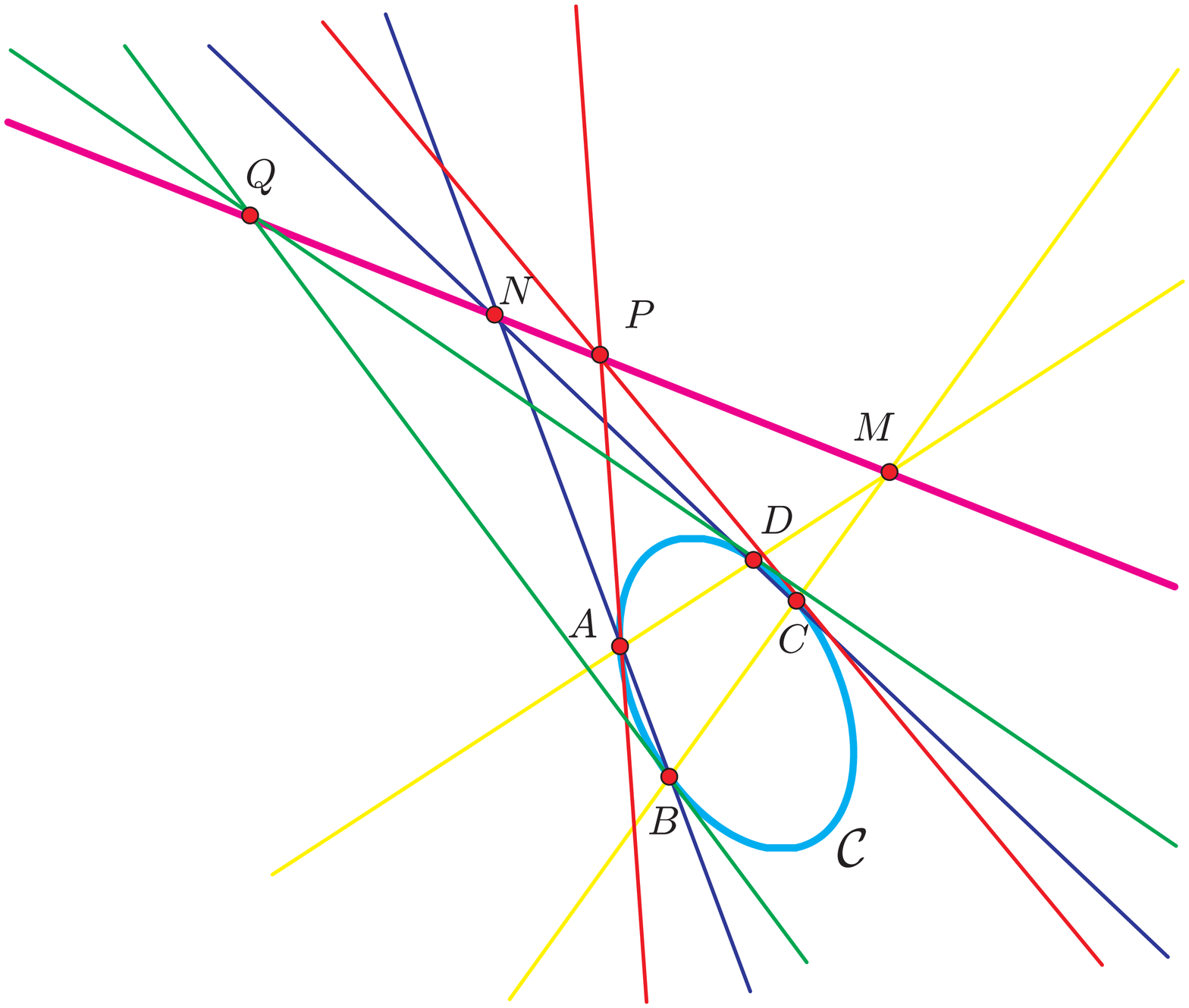,
width=0.6\textwidth}} \caption{Proposition \ref{slikacetvorougao}}
\label{cetvorougaouelipsi}
\end{figure}

\begin{proposition}\label{slikacetvorougao} Let $A B C D$ be a quadrilateral inscribed in a conic $\mathcal{C}$ and let the point $M$ be the intersection of the lines $A D$ and $B C$, the point $N$
the intersection of the lines $A B$ and $C D$, the point $P$ the
intersection of the tangents to $\mathcal{C}$ at $A$ and $C$, and
$Q$ the intersection of the tangents to $\mathcal{C}$ at $B$ and
$D$. Then the points $M$, $N$, $P$ and $Q$ are collinear (see
Figure \ref{cetvorougaouelipsi}).
\end{proposition}

\begin{figure}[h!h!]
\centerline{\epsfig{figure=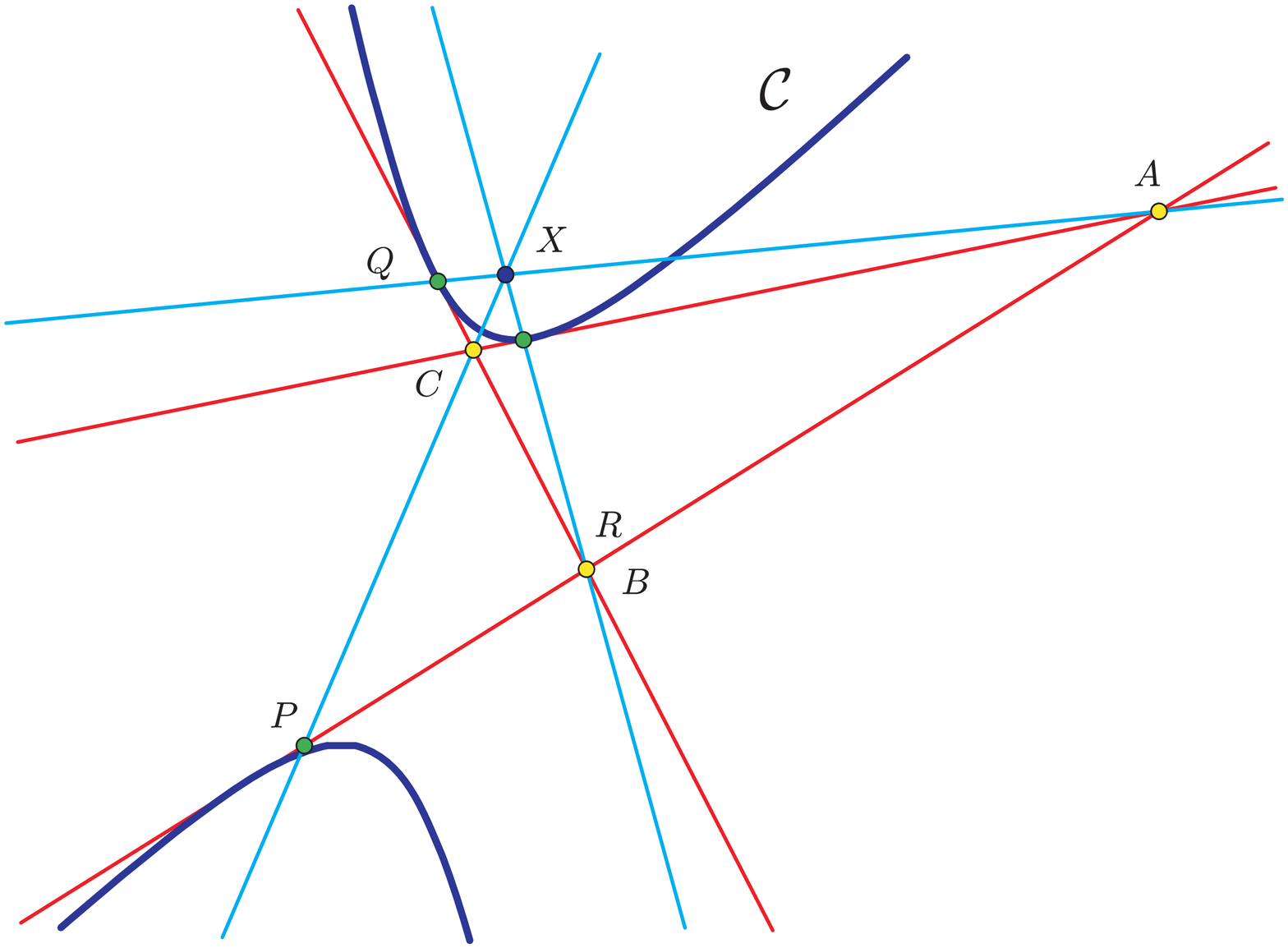,width=0.7\textwidth}}
\caption{Proposition \ref{tro}} \label{trougao}
\end{figure}

\begin{proposition}\label{tro} Let $A$, $B$ and $C$ be points such that the lines
$A B$, $B C$ and $C A$ tangent conic $\mathcal{C}$ in the points
$P$, $Q$ and $R$, respectively. Then the lines  $A Q$, $B R$ and
$C P$ belong to the same pencil of lines (see Figure
\ref{trougao}).
\end{proposition}

There are many ways of obtaining a degenerate configuration. If we
take in classical octagrammum mysticum $G\rightarrow E$ and
$H\rightarrow F$ then it is possible to get classical Pascal
theorem. Thus far, the degeneracy tool is a machinery for getting
many geometrical theorems about $n$-gons inscribed in a conic.

\begin{figure}[h!h!]
\centerline{\epsfig{figure=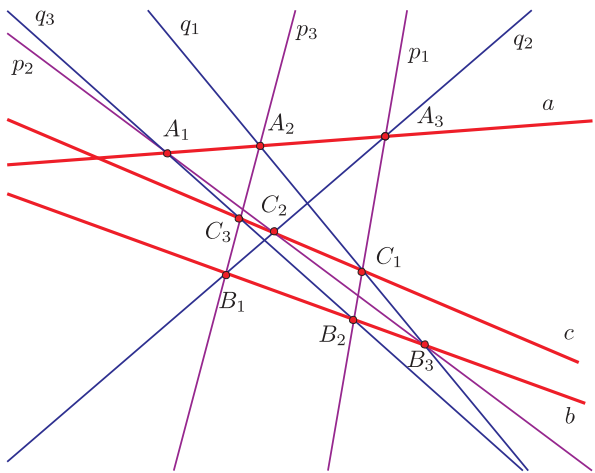,width=0.6\textwidth}}
\caption{Pappus theorem} \label{papos}
\end{figure}

Also it is possible to take conic $\mathcal{C}$ to be degenerate.
Then we obtain the theorem of Pappus and its generalizations (see
Figure \ref{papos}). Pappus's theorem was the first result about
mystic hexagon. Applying this theorem several times it is possible
to obtain dynamical system which is described completely in
\cite{Schwa}.

\section{Connections with Other Constructions in Geometry and
Combinatorics}

In the previous sections we focused on proving some interesting
geometrical facts about hexagon and octagon inscribed in a conic.
Our main technic was essentially a careful application of
Corollary \ref{vazno}. Now we want to put these configurations in
some other context.

Paper \cite{Ladd} gives a beautiful combinatorial description of
hexagons generating Kirkman, Steiner and Salmon points, as well as
Kirkman and Steiner lines. From this it is easy to conclude which
subgroups of $S_6$ are associated with the lines and the points in
configurations. In \cite{Carde} it is explained how outer
automorphisms of the permutation group $S_6$ reflect on the
Pascal's configuration.

It is natural to ask what's happen in case of octagon inscribed in
a conic. As seen, we could obtain in general case 2520 conics in
configuration. Start with the conic $\mathcal{C}^\ast$ arising
from the quartics $l (A B)\cdot l (C D) \cdot l (E F) \cdot l (G
H)$ and $l (B C) \cdot l (D E) \cdot l (F G) \cdot l (H A)$. The
first natural question is which subgroup of $S_8$ fixes this conic
with respect to an action of $S_8$ on the vertices of octagon.

\begin{proposition} The mystic conic $\mathcal{C}^\ast$ is fixed with the dihedral
subgroup $D_8$ of $S_8$ which generators are $8$-cycle and
reversing order permutation.
\end{proposition}

It is also natural to ask what happen with conics arising from two
quadrilaterals inscribed in a conic. Such conic is stable with
respect to $D_4\times D_4$ subgroup of $S_8$. If we are talking
about subgroups fixing some fixed pencil of conics in octagrammum
mystic, situation is complicated. We distinct pencils by groups
fixing it.

\begin{proposition}\label{tip1} The pencil of conics in octagrammum
mysticum formed with the quartics $l (A B)\cdot l (C D) \cdot l (E
F) \cdot l (G H)$, $l (B C) \cdot l (D E) \cdot l (F G) \cdot l (H
A)$ and $l (A D) \cdot l (C H) \cdot l (G E) \cdot l (F B)$ is
fixed with certain subgroup $Z_3 \times D_8$ of $S_8$.
\end{proposition}

\begin{proposition}\label{tip2} The pencil of conics in octagrammum
mysticum formed with the quartics $l (A B)\cdot l (C D) \cdot l (E
F) \cdot l (G H)$, $l (B C) \cdot l (D E) \cdot l (F G) \cdot l (H
A)$ and $l (A F) \cdot l (C E) \cdot l (D G) \cdot l (B H)$ is
fixed with subgroup $D_8$ of $S_8$ that fixes all three conics.
\end{proposition}

In fact Propositions \ref{tip1} and \ref{tip2} clearly explain the
difference between pencils. On each conic from mystic octagon
there are exactly $2$ conics that are fixed by some subgroup $Z_3
\times D_8$ of $S_8$. So in fact there are only 1680 such pencils.
Other pencils are fixed with subgroup $D_8$ of $S_8$ that fixes
every conic.

Now we will switch attention to the other field. In \cite{Yuz} the
notion of $(k, d)$ nets is defined. Yuzvinsky proved that only
possible values for $(k, d)$ are: $(k=3, d\geq 2)$, $(k=4, d\geq
3)$, $(k=5, d\geq 6)$. Examples and constructions of some $3$-nets
are given in \cite{Stip}. Let observe that in hexagrammum mystic
there is also example of $(3, 4)$ net. Take in Theorem \ref{gsl}
$\mathcal{D}=l (A F)\cdot l (B E)\cdot l (C D)$. Sets of lines
$\{p_1, p_2, p_3, p_4\}$, $\{q_1, q_2, q_3, q_4\}$ and $\{l (A F),
l (B E), l (C D), s\}$ where $s$ is Steiner line.

Let us mention that the theory of algebraic curves is closely
related with problems treating $(k, d)$-nets of lines. From
theorem \ref{gslc} we will construct one examples with conics.

\begin{example} Take in Theorem \ref{gslc} for quartic
$\mathcal{Q}=\mathcal{C}_4\cdot \mathcal{D}_4$ where
$\mathcal{C}_4$ is a conic through points $A$, $B$, $C$ and $D$
and $\mathcal{D}_4$ is a conic through points $E$, $F$, $G$ and
$H$. Consider the following sets of conics $\{\mathcal{X}_1,
\mathcal{X}_2, \mathcal{X}_3\}$, $\{\mathcal{Y}_1, \mathcal{Y}_2,
\mathcal{Y}_3\}$ and $\{\mathcal{C}_4, \mathcal{D}_4,
\mathcal{S}\}$ where $\mathcal{S}$ is the conic obtained by
Theorem \ref{gslc}. Then this set are examples for $3$-net of
conics, which precise definition we will give analogously with the
definition of $k$-nets of lines.
\end{example}

\begin{definition}\label{defcon} Let $k$ be a positive integer and $P$
projective plane. A $k$-net of conics in $P$ is a $(k+1)$-tuple
$(\mathcal{A}_1, \dots, \mathcal{A}_k, \mathcal{X})$, where each
$\mathcal{A}_i$ is a nonempty finite set of conics in $P$ and
$\mathcal{X}$ finite set of pencils of conics, satisfying the
following conditions:
\begin{enumerate}
    \item The $\mathcal{A}_i$ are pairwise disjoint.
    \item If $i\neq j$, then the pencil of conics generated by any
    conic from $\mathcal{A}_i$ and any conic from $\mathcal{A}_j$
    belongs to $\mathcal{X}$.
    \item Trough every pencil in $\mathcal{X}$ passes exactly one
    conic of each $\mathcal{A}_i$.
\end{enumerate}
\end{definition}

This definition is uninteresting when $k=2$. Any two disjoint sets
of conics form $2$-net. When $k\geq 3$ then as in case of $k$-net
for line, following same idea, see for example \cite{Stip} we
obtain the same combinatorial restriction for $k$-net of conics:

\begin{proposition}\label{moj1} Let $(\mathcal{A}_1, \dots, \mathcal{A}_k,
\mathcal{X})$ be a $k$-net of conics, with $k\geq 3$. Then every
$\mathcal{A}_i$ has the same cardinality. Furthermore, pencils of
$\mathcal{X}$ are generated by conics of $\mathcal{A}_i$ and
$\mathcal{A}_j$, for any $i\neq j$. Thus
$|\mathcal{X}|=|\mathcal{A}_i|^2$.
\end{proposition}

Following terminology for the case of lines, the cardinality of
$\mathcal{A}_i$ is called \textit{degree} of the $k$-net of
conics.

Note that definition \ref{defcon} could be reformulated in the
following sense. Since $P^5$ is the space of conics in projective
plane, and pencil of conics are the lines in $P^5$, then we could
look at $\mathcal{A}_i$ as sets of points in $P^5$ and
$\mathcal{X}$ as set of lines in $P^5$ such that:
\begin{enumerate}
    \item The $\mathcal{A}_i$ are pairwise disjoint.
    \item If $i\neq j$, line through any point from $\mathcal{A}_i$ and any from $\mathcal{A}_j$
    belongs to $\mathcal{X}$.
    \item On every line in $\mathcal{X}$ lies exactly one
    point of each $\mathcal{A}_i$.
\end{enumerate}

One of thinking about the Pascal's hexagon and octagon
configurations is to treat them like arrangements of curves in
$\mathbb{C} P^2$. This view is quite present in contemporary
research on the subject and it would be interesting to apply some
results concerning the arrangements of curves to the case of
mystic hexagon and octagon configurations. We believe the
invariants like the Solomon-Orlik algebra and the cohomology of
the complement of arrangements could give interesting results.

\section{Further Research}

Here we will discus some questions that could be of interest for
the further research.

\textbf{Question 1.} \textit{Is there any interesting conic that
belong to several pencil of conics in the octagrammum mystic
except one found in Theorem \ref{gslc}?}

As one can see from the proof for Salmon-Cayley line this is not
question that could be easily answered. In fact if we don't have
natural candidate it is hard to smartly apply B\'{e}zout's
theorem. Maybe the difference between some pencils established in
Propositions \ref{tip1} and \ref{tip2} give some hope to the
affirmative answer due to an analogy with hexagon case.

\textbf{Question 2.} \textit{What the theory of arrangements of
curves says about the mystic hexagon and octagon?}

In this moment, this is just idea how to look on problem in
context of modern mathematics. We don't have any assumption what
in fact we expect from higher technics and theories. But
nevertheless, we believe that beautiful theorem about hexagon and
octagon inscribed in conic could be revealed in new unexpected
shape.

\textbf{Question 3.} \textit{Determine all the possible values $k$
and $d$ such that $(k, d)$-net of conics exist?}

At first, answer for the net of lines is not given completely. But
on the first look much of things that are done for case of lines
could be tried in the case of conics. This is problem we will do
in future.

\textbf{Question 4.} \textit{Find the new examples of $k$-nets of
conics.}

In fact we only gave one example of $(3, 3)$-net of conics.
Construction of $(3, d)$ net of lines is strongly connected with
orthogonal Latin squares and has interesting combinatorial
structure. Thus finding new examples in the case of conics and
some method for generating such examples will be interesting.

\textbf{Question 5.} \textit{Describe the dynamical system
obtained by the application of the octagon mysticum.}

This is the question we posed having in mind the famous result of
Schwartz, see \cite{Schwa}, where one particular dynamical system
arising from Pappus theorem is explained by modular group. The
system that naturally arises from an octagon inscribed in a conic
is much richer because if we treat Theorem \ref{8conc} like move
then it could be implied in many ways. We hope to give answer
these questions in the near future.

\textbf{Question 6.} \textit{When the octagrammum mysticum
produces ellipse, hyperbola, parabola and when degenerate conics?}

We worked in the space $\mathbb{C} P^2$. From the standard
Euclidean picture this is natural and hard question. But something
it is done in paper \cite{Gold}, so we hope that is possible to
make some advance toward this question.

\begin{center}\textmd{Acknowledgements }
\end{center}

The authors are grateful to the participants of Belgrade seminars:
Mathematical Methods of Mechanics, Geometry Seminar and CGTA
Seminar for useful suggestions and discussions, as well as their
constant support.

\bigskip

{\small \DJ{}OR\DJ{}E BARALI\'{C}, Mathematical Institute SASA,
Kneza Mihaila 36, p.p.\ 367, 11001 Belgrade, Serbia

E-mail address: djbaralic@mi.sanu.ac.rs

\bigskip
IGOR SPASOJEVI\'{C}, University of Cambridge, Cambridge, United
Kingdom

E-mail address: igoraspasojevic@yahoo.co.uk}


\begin{thebibliography}{99}

\bibitem{Bach} I. Bacharach, \textit{\"{U}ber den Cayley'schen Schrittpunktsatz}, Math.
Ann. 26 (1886), 275-299.

\bibitem{Ber}  M. Berger, \textit{Geometry
Revealed - A Jacob's Ladder to Modern Higher Geometry}, Springer
Heidelberg Dodrecht London New York, 2010.

\bibitem{Bix} R. Bix, \textit{Conics and Cubics - A Concrete Introduction
to Algebraic Curves}, Springer, 2006.

\bibitem{Boro} M. Borodzik, H. \.{Z}o\l \c{a}dek, \textit{The Pascal
Theorem and Some Its Generalizations}, Topological Methods in
Nonlinear Analysis, Volume 19, 2002, 77-90.

\bibitem{Bott} O. Bottema, \textit{A Generalization of Pascal's theorem,
Duke Mathematical Journal}, Volume 22, 1955, 123-127.

\bibitem{Carde} H. C\'{a}rdenas, R. San Agust\'{i}n, \textit{On Veronese's
Decomposition Theorems and the Geometry of outer automorhisms of
gropus $S_6$}, Journal of Combinatorial Mathematics and
Combinatorial Computation, 1997

\bibitem{Carv} W. B. Carver, \textit{On the Cayley-Veronese Class of
Configurations}, Transactions of the American Mathematical
Society, Volume 6, No. 4, 1905, 534-545.

\bibitem{Cayl} A. Cayley, J. f. Reine u. Angewandte Math. 31. 1864, 213.

\bibitem{Cay} A. Cayley, \textit{On the Intersection of Curves}, Cambridge Math. J. 3
(1843), 211-213.

\bibitem{Clark} B. G. Clark, \textit{An Analytic Study of the Pascal
Hexagon}, American Mathematical Monthly, 1937, 228-231.

\bibitem{Cour} N. A. Court, \textit{Pascal's Theorem in Space}, Volume 20, 1953, 417-421.

\bibitem{Fox} C. Fox, \textit{The Pascal Line and its Generalizations}, American Mathematical
Monthly, Volume 65, No. 3, 1958, 185-190.

\bibitem{Grif} P. A. Griffiths, \textit{Introduction to Algebraic
Curves}, American Mathematical Society, 1989.

\bibitem{gh} P. Griffiths, J. Harris, \textit{Principles of Algebraic
Geometry}, John Wiley \& Sons, Inc., 1978.

\bibitem{Katz} G. Katz, \textit{Curves in Cages: an Algebro-geometric Zoo}, Amer.
Math. Monthly, 113(9), 2006, 777-791.

\bibitem{Gold} L. S. Evans, J. F. Rigby,
\textit{Octagrammum Mysticum and the Golden Cross-Ratio}, The
Mathematical Gazette, Vol. 86, No. 505 (Mar., 2002), pp. 35-43.

\bibitem{Kirkman} T. P. Kirkman, Cambridge and Dublin Journal, Volume V 1850,
185.

\bibitem{Kir} F. Kirwan, \textit{Complex Algebraic Curves}, Cambridge
University Press, 1992.

\bibitem{Ladd} C. Ladd, \textit{The Pascal Hexagram}, American Journal of
Mathematics, Volume 2, 1879, 1-12.

\bibitem{Pasc} B. Pascal, \textit{Essai por les coniques},  Eubres, Brunschvigg et
Boutroux, Paris, Volume I, 1908, 245.

\bibitem{Pra}  V. V. Prasolov and V. M. Tikhomirov,
\textit{Geometry}, American Mathematical Society, 2001.

\bibitem{Gerb}  J. Richter-Gebert,
\textit{Perspectives on Projective Geometry}, Springer-Verlag
Berlin Heidelberg, 2011.

\bibitem{Rupp} C. A. Rupp, \textit{An Extension of Pascal's Theorem},
Transactions of the American Mathematical Society, Volume 31, No.
3., 1929, 580-594.

\bibitem{Sal} G. Salmon, \textit{A Treatise on the Analytic Geometry of
Three Dimensions}. 5th ed., New York, 1912.

\bibitem{Salm} G. Salmon, \textit{A treatise on Conic Sections}
Chelsea, 1952, 219-236.

\bibitem{Schwa} R. Schwartz, \textit{Pappus's theorem and the
modular group}, Publications math\'{e}matiques de l\'{e}Institut
des hautes \'{e}tudes scientifiques, 78, 1993, 187-206

\bibitem{Stei} J. Steiner, Gesamelte Werke, Volume I, 1832, 451.

\bibitem{Stip} J. Stipins, \textit{OLD AND NEW EXAMPLES OF $k$-NETS IN
$\mathbb{P}^2$}, arXiv:math/0701046v1

\bibitem{Thom} H. D. Thompson, \textit{A Note on Pencils of Conics},
American Journal of Mathematics, Volume 9, No. 2, 1887, 185-188.

\bibitem{Vero} G. Veronese, \textit{Nuovi Teoremi sull' Hexagrammum Misticum},
Memorie dei Reale Accademie dei Lincei, Volume I, 1877, 649-703.

\bibitem{Wilk} T. T. Wilkinson, Mathematical questions with their solutions,
2015, Educational Times XVII (1872) p.72.

\bibitem{Yuz} S. Yuzvinsky, \textit{Realization of finite abelian groups
by nets in $\mathbb{P}^2$}, Compos. Math, 140, no. 6, 2004
1614-1624

\end{thebibliography}
\end{document}